\documentclass[12pt]{amsart}
\usepackage{amsmath}
\usepackage{listings}
\usepackage{amsfonts,amssymb,amsthm, txfonts, pxfonts}
\usepackage[title,titletoc,toc]{appendix}
\usepackage{lipsum}
\newcounter{quotecount}

\begin{document}

\title[Igor Rivin]{Some Thoughts on the Teaching of Mathematics -- ten years later}
\author{Igor Rivin}
\address{Department of Mathematics, Temple University, Philadelphia, PA}
\curraddr{Mathematics Department, Brown University, Providence, RI}
\email{rivin@temple.edu}
\email{igor\_rivin@brown.edu}
\thanks{The author was supported in part by the NSF DMS when the original version was written, and is now supported by ICERM under their block NSF grant. These thoughts are
  very much in a preliminary   DRAFT form. The author would like the many people who commented on the previous version of this document.I would like to thank the anonymous referees for helpful suggestions} 
\date{\today}
\maketitle

%\section{Introduction}
This is a somewhat expanded and corrected version of a ``manifesto'' first posted on my Temple web page almost exactly 10 years ago. Since then, I have put some of the ideas expressed below into practice, -- a brief description of my experience is included in Section \ref{experience},  and I have also had a fair amount of feedback. Some (most) of it has been completely positive, others have included some critical ideas. I describe some of it in Section \ref{feedback}.

\section{What is the problem?}

A mathematics professor in a public university has many
responsibilities. These include research, administration, and
teaching. Teaching, in turn, include ``specialized'' teaching (to mathematics
majors and graduate students) and ``service'' teaching: teaching mathematics
to first and second year students. These thoughts will center primarily on
service teaching, which, for me, combines some of the most exciting and some
of the most depressing aspects of my job. Some thoughts on teaching mathematics majors are added in Section \ref{specialist}.

\subsection{The product} Why depressing? Consider: the vast majority of
service courses are concerned 
with differential and integral calculus and linear algebra. These are both
rather deep subjects, as evidenced by the fact that mathematics had been
practiced for thousands of years by rather talented people before the basic
principles of the calclulus were laid down in the late seventeenth century (although some of these principles were discovered, in an \emph{ad hoc} way by Archimedes, considerably earlier). It
took another two hundred years of extensive work to make the foundations of
the subject truly solid. Linear algebra, as used today, is an even later
bloomer. The current machinery of matrices and linear transformations was not
put into a truly modern form until the beginning of the twentieth century.  We have no choice but to agree
hat these subjects are quite deep and require
some considerable technical skill to use successfully.

\subsection{The consumers} \emph{Who} are we teaching them to? In a public university (such as Temple) our
students are, in the main, somewhat above-average products of the US public
secondary education system. This mean that their techincal ability is already quite severely taxed by arithmetic with fractions.
Their abstract reasoning skills are essentially nonexistent, and the very
concept of proof is foreign to them. 

\subsection{The results} A consequence of all this is that it is well-nigh \emph{impossible} to teach
them what 
we purport to teach them: higher mathematics presupposes a certain level of
abstraction, and even if we commit the crime of forgetting that, and define
calculus as ``a collection of computational techniques without
understanding,'' the students' technical weakness renders even that aspect
essentially worthless. They \emph{can not} compute. The result is that our
calculus and linear algebra classes consist of a collection of trivial
examples which the students must memorize by rote. This has the consequence
of not teaching the students anything except the fear and hatred of
mathematics. There is more still: the majority of the students never use
calculus in their future lives (small wonder, since they don't actually know
any, as discussed above,) but they never had any intention of using advanced
mathematics even \emph{before} taking the courses). They are required to take
the courses because of the (not unreasonable) belief that mathematics should
be a part of ever college-educated person's intellectual make-up. The result
is that the \emph{loathing} of mathematics is part of the intellectual makeup
of a sizeable majority of Americans. The amazing (and exhilarating)
observation is that despite all of the above, some students actually manage
to understand something of the subject. The exhilaration is, however,
tempered by the thought of the huge amount of wasted time and by the thoughts
of what these talented students could achieve if taught properly.\footnote{A typical cafe flirtation usually starts with: \textsc{Girl:} what do you do? \textsc{Mathematician:} I study mathematics. \textsc{Girl:} (one of the two responses: I am terrible at math OR I used to love math!). This is a very American phenomenon; I have never had such a conversation in Europe.}

\section{What can we do?}
What, then, is the solution? We could drop the distribution requirements in
Mathematics (and I could easily see this happening), but the fact of the
matter is that the ability to reason logically and abstractly really should be
something (perhaps the \emph{main} thing) everyone takes from their higher education, and something we, as mathematics
educators (which we are, even if the term does produce a visceral reaction in
most people),  ought to instill it in our students. How? The first step is a step back
-- a step back from ``higher'' mathematics -- the mathematics of infinite and
the infinitesimal --  back to conceptually simpler forms of
mathematics. 

In the (not so distant) past, Euclidean geometry was such a subject, taught
exactly for the above-stated reasons (the fundamental concepts -- of line,
circle, distance, \emph{etc.} are quite intuitive, while the basic components
of mathematical reasoning are all present). I would not necessarily advocate a return to
this, however.  Firstly, the subject has been dead for several hundred
years and secondly, it is quite far from the modern American experience.

Instead, it is my opinion that we should start at the very beginning -- with
reasoning (logic) and counting (which means naive set theory and combinatorics
and graph theory) and probability. These subjects are ever more visibly
important in our lives due to the ubiquity of computers. They are
both easier to learn and more immediately rewarding for the students than what we are currently teaching. In
addition, there is another principle which can be used to clear at least some
of the mush out of the students' minds. That principle is: 
\begin{equation}
\label{programproof}
\boxed{
\textsc{\mbox{A computer
program \emph{is} a proof.}}}
\end{equation}
For example, the student who wrote the program in Section \ref{crt} had to have complete understanding of the Chinese Remainder Theorem, and his program \emph{is} the chinese remainder theorem, in that, given some quantities satisfying the hypotheses of the theorem, it never fails to produce a quantity satisfying the conclusion.

This correspondence between program and proof (though far from perfect)
allows us to make mathematics hands on (programming a computer gives very
rapid feedback, both positive and negative), it is closely related to the students'
experience and visibly ``useful.'' (of course, the real utility of mathematics
lies much deeper, but \dots).

One problem with this approach is the potential need to waste a lot of time
introducing students to the subtleties and idiosyncrasies of some (possibly
proprietary) programming language or a scientific computing system. It is
important to start the students off on a mathematically clean system,
preferably running on a mathematical  \emph{abstract   machine}. Luckily, such
a system is available and had been used for over thirty years in
\emph{Computer Science} education at MIT, University of Indiana, and many
other schools. This is the SCHEME programming language.\footnote{See
  \texttt{http://www.scheme.org} for more details.} The justly acclaimed
book by H.~Abelson and G.~J.~Sussman (``Structure and Interpretation of
Computer Programs'', \cite{abelsonstructure}) introduces the fundamentals of computer programming, and
together with a companion book on ``Structure and Interpretation of
Mathematics''\footnote{Writing such a book is a project I am very keen on, but
  this has to be done in parallel with using these ideas in teaching -- the
  book of Abelson and Sussman was circulated as lecture notes for several
  years}
this would constitute the core of a modern introduction to
Mathematics.  It should be noted that SCHEME is no longer used as the introduction to computing at MIT, having been supplanted by Python, but the reason for this is that Python is a more useful tool for practical work (having libraries for almost any task one might need to accomplish), and has absorbed much of the SCHEME semantics. It is, however, worse as a first language, because it is less pure, and has less well-defined semantics.

Given the very high level of functionality presented by \texttt{Python,} \texttt{Mathematica,} or other "kitchen sink" languages, it may be a reasonable compromise to start with Scheme as a way of building a foundation and then proceeding to use one of these languages when (more precisely, when and if) more advanced topics mathematics needs to be introduced.

\section{Possible Objections}
One could foresee some objections to the above program:

\subsection{But what about calculus?} There is no suggestion of eliminating
calculus completely from the University curriculum: the students of the
sciences and engineering (not to mention mathematics) do need to be acquainted
with it. However, it would enter somewhat later, when the students are more
mathematically mature, and thus more capable of actually understanding
something of the subject. It is true that a lot fewer students will be
\emph{required} to study calculus, but, on the other hand, the number of
people studying \emph{Mathematics} may actually increase.

\subsection{But what about the teachers?} It has been my observation that many
of the TAs, having been educated in a ``traditional'' way are a little rusty
on ``finite'' mathematics. Same goes (in spades) for many of the
professors. Many faculty members and graduate students have no familiarity
with computer programming at all. The first problem is easily fixable, the
second slightly less so, but my contention is that anything which is teachable
to freshmen should be even easier to teach to faculty (and is no less useful
to them).

\section{My experience over the last decade}
\label{experience}
Since the first version of this document was published, I taught a number of relevant courses (with multiplicity). One was ``Junior Problem Solving'', which is an introductory course in mathematics for computer science majors. Another was ``Mathematical patterns'', which is a sort of a ``Math for Poets'' course, a third was ``Senior problem solving'', which is (in principle) the last mathematics course a mathematics major at Temple takes. The fourth was ``mathematical computing'' -- a course for  mathematics graduate students. 
\subsection{Poets}
In the first two courses (Junior Problem Solving and Mathematical Patterns) I taught the basics of logic, as presented in the very nice book \cite{gensler2010introduction} by Harry Gensler ,who is an ethical philosopher, not a mathematician. This means, in particular, that while the book eventually goes into symbolic logic, even there, a lot of the questions are about deciding validity of English sentences or (more generally) philosophical arguments. In both instances, the course was a success, in that, while, as far as I could tell, \emph{none} of the students  could carry out a logical argument at the start of the course, \emph{many} could at the end (the fact that they could not at the beginning validates my belief that the students were not ready for anything resembling calculus).
Some additional remarks:
\begin{itemize}
\item At one point in the Junior Problem Solving course, I decided to try something more "mathematical", and presented the proof that $\sqrt{2}$ is irrational. Despite going very slowly, it was quite clear that the arithmetic involved was too much for the students.
\item Many of the students thanked me at the end of the course (in one case, I was walking down the street in Center City Philadelphia, when one of my former students crossed the street, dodging a number fast moving cars, to tell me how much she appreciated the material). This sort of thing had never happened at the conclusion of any calculus course I had taught.
\end{itemize}
\subsection{Budding Mathematicians}
\label{specialist}
In the graduate course, I started by teaching the elements of programming \emph{a la} Abelson-Sussman, with more emphasis on mathematical problems. The experience was definitely bimodal. The good students (by which I mean, ones having a lot of mathematical potential) had extremely good programming skills, while the not-so-good students had essentially none. This, of course, is quite unfortunate, since most of Temple's PhD candidates would wind up in the industry, which means that programming skills are essential to their future livelihood -- the best students are the ones most likely to pursue a career in pure mathematics, and in that sense \emph{need} programming skills the least (although , as Gauss already knew well, mathematics is an experimental science, and computing is the experiment).

The Senior Problem Solving course was, in a way, the most disappointing -- while some of the students were extremely good, most did not have much better technical skills than the computer science students in Junior Problem Solving (to make it more depressing, these were mostly mathematical education majors), and would tend to get just as confused by proofs, especially ones which required computation. I am quite sure that had they been required to take a course of the kind I was teaching to their poetical brethren, together with a rigorous algebra course (of the sort those of us who grew up in the Soviet Union took at the age of 12) they would have been much more deserving of a mathematics degree. However, by the time I got them, it was already too late.
%specialist section should be a subsection of this section.
%\begin{appendices}
\appendix
\section{Source code for a modular arithmetic package in \texttt{Scheme}}
\label{crt}
The code below (taken verbatim from a homework assignment done by a student in one of my classes) defines all of modular arithmetic in Scheme (which already has most of it, so such a program would be a lot shorter and more efficient in practice). At the end, an implementation of the extended Euclidean algorithm is given, followed by an implementation of the Chinese remainder theorem algorithm.
\lstset{language=Lisp} 
\begin{lstlisting}
(define (intMod k m)
  (if (= m 0)
      (lambda (s) 0)
      (let* ((j (remainder k m))
             (l (if (< j 0)
                    (+ j m)
                    j)))
        (lambda (s)
          (if s 
              l
              m)))))
  
(define (sameMod? x y)
  (= (x #f) (y #f)))

(define (eqM? x y)
  (and (= (x #f) (y #f))
       (= (x #t) (y #t))))

(define (displayM x)
  (display (x #t)) 
  (display " mod ") 
  (display (x #f)) 
  (display "\n"))
          
(define (+M x y)
  (if (sameMod? x y)
      (intMod (+ (x #t) (y #t)) (x #f))
      (intMod 0 0)))

(define (-M x y)
  (if (sameMod? x y)
      (intMod (- (x #t) (y #t)) (x #f))
      (intMod 0 0)))

(define (*M x y)
  (if (sameMod? x y)
      (intMod (* (x #t) (y #t)) (x #f))
      (intMod 0 0)))


(define (Modulus x)
  (x #f))

(define (Coset x)
  (x #t))

(define (gcde a b) 
  (letrec ((aux (lambda (x1 x2 y1 y2 r1 r2)
                  (let ((r3 (remainder r1 r2)) (q (quotient r1 r2)))
                    (if (= r3 0)
                        (list r2 x2 y2)
                        (aux x2
                             (- x1 (* q x2))
                             y2
                             (- y1 (* q y2))
                             r2
                             r3))))))
    (aux 1 0 0 1 a b)))

(define (/M x y)
  (let ((z (apply gcde (list (y #t) (y #f)))))
  (if (and (sameMod? x y) (= 0 (remainder (x #t) (car z)))) 
      (intMod (* (cadr z) (/ (x #t) (car z))) (x #f))
      (intMod 0 0))))

  
(define (crt l)
  (if (null? (cdr l))
      (car l)
      (let* ((w (map Coset (list (cadr l) (car l))))
             (x (map Modulus (list (car l) (cadr l)))) 
             (y (apply gcde x))
             (z (car y)))
        (if (apply = (map remainder w (list z z)))
            (crt (cons 
                  (intMod (/ (apply + (map * w x (cdr y))) z) 
                          (/ (apply * x) z)) 
                  (cddr l)))
            (intMod 0 0)))))
\end{lstlisting}
\section{Some feedback}
\label{feedback}
Most feedback has been overwhelmingly positive. Here are some examples:
\begin{itemize}
\item Some have suggested probability as a good introduction to discrete math. Ward  and Gundlach (of Purdue) have written a book for the purpose.
\item
\begin{quotation}
I read your diatribe and found it very compelling---in fact, something that I agree with and find in accordance with my own experience as an undergraduate (an electrical engineering major). It wasn't until I took a discrete mathematics course in my Sophomore year that I began to understand mathematics. I'm now at Caltech studying applied math, so to say that discrete math was the start of an enormous change in my career path would not be an understatement.

At the same time, you will experience a large number of engineers who feel that a lack of, say, vector calculus in 3 dimensions will seriously hinder students in (say) electromagnetic engineering. I found this to be the case even though I had taken the course and ostensibly "knew" the material. 
\end{quotation}
Of course, I certainly agree that engineers should study multivariate calculus, but not as the first thing.
\item
It was suggested that the ``Integrated approach'' (where mathematics and the science that uses it are used in parallel) is a solution to some of the problems I am trying to address. My view on this is that while this is an excellent idea, it cannot also be done in parallel to teaching people to think.
\item (from a chemist at Temple):
\begin{quotation}
Hey, I was reading the first part of your diatribe on mathematics education and so far I agree with everything you said. I feel the same way about teaching Organic Chemistry - I feel like the course is more about identifying the 5% of the population that can handle Organic Chemistry than teaching.
\end{quotation}
\item (from an applied mathematician friend):
So I ask you:
\begin{quotation}
        What is the value of an education system which starts at "the
        very beginning" at 18 years old ?
\end{quotation}

My answer is: very little.

So.  Can we do better?  The answer is *yes*.  Try google on
"Montessori".  Both my children spent from age 3 to 12 in Montessori
schools.  When I try to explain any logical arguments to them they look
at me as if I am a complete moron because by age 12 in a Montessori
system, logical reasoning is natural.
\end{itemize}
%\end{appendices}

\bibliographystyle{plain}
\bibliography{teachingnew}
\end{document}